\documentclass[10pt,twocolumn]{article}
\usepackage{algorithm}
\usepackage{hyperref}
\usepackage[noend]{algpseudocode}
\makeatletter
\def\BState{\State\hskip-\ALG@thistlm}
\makeatother
\usepackage{amsmath}
\usepackage{amssymb}
\usepackage{graphics}
\usepackage{graphicx}
\usepackage{oldlfont}
 \usepackage{newlfont}
 \usepackage{amsfonts}
 \usepackage{amssymb}
 \usepackage{mathtools}
\usepackage[normalem]{ulem}
\usepackage{color}
 \usepackage{ulem}
 \newcommand{\R}{{\mathbb R}}
   \usepackage{amsthm}
 \theoremstyle{definition}

\usepackage[table]{xcolor}
\usepackage{mathrsfs}
\usepackage{multirow}
\newcommand{\ie}[0]{\textit{i.e.}}

\newcommand{\mbArgmin}{\operatornamewithlimits{argmin}}

\newcommand{\mbA}{\mathbf{A}}  % Measurements matrix
\newcommand{\mbW}{\mathbf{W}} % Shearlets matrix
\DeclareMathOperator{\sign}{sign}

\newcommand{\Z}{\ensuremath{\mathbb{Z}}}

\def\lra{\longrightarrow}
\def\<{\langle}
\def\>{\rangle}

\newcommand{\ve}[1]{\boldsymbol{#1}}
\def\f{{\ve f}} % object
\def\m{{\ve m}} % measurements/sinogram
\def\y{{\ve y}} % dual
\def\v{{\ve v}} % dual

\begin{document}

\title{Controlled Wavelet Domain Sparsity for X-ray Tomography}

\author{Zenith~Purisha$^{1,2}$, Juho~Rimpel\"{a}inen$^1$, Tatiana Bubba$^1$ and~Samuli~Siltanen$^1$}
\date{%
    $^1$Department of Mathematics and Statistics, University of Helsinki, Finland\\%
    $^2$Department of Mathematics, Universitas Gadjah Mada, Indonesia\\[2ex]%
    \today
}
\maketitle

%\tableofcontents
\begin{abstract}
Tomographic reconstruction is an ill-posed inverse problem that calls for regularization. One possibility is to require sparsity of the unknown in an orthonormal wavelet basis. This, in turn, can be achieved by variational regularization, where the penalty term is the sum of the absolute values of the wavelet coefficients. The primal-dual fixed point (PDFP) algorithm introduced by Peijun Chen, Jianguo Huang, and Xiaoqun Zhang (Fixed Point Theory and Applications 2016) showed that the minimizer of the variational regularization functional can be computed iteratively using a 
soft-thresholding operation. Choosing the soft-thresholding parameter $\mu>0$ is analogous to the notoriously difficult problem of picking the optimal regularization parameter in Tikhonov regularization. Here, a novel automatic method is introduced for choosing $\mu$, based on a control algorithm driving the sparsity of the reconstruction to an {\it a priori} known ratio of nonzero versus zero wavelet coefficients in the unknown. 
\end{abstract}
\\
\noindent
{\bf {\it Keywords}} : tomography, wavelet, sparsity, regularization, control, limited data tomography, X-ray

% \begin{keywords}
% computed tomography, wavelet, sparsity, regularization, automatic controlled
% \end{keywords}

%\IEEEpeerreviewmaketitle

\section{Introduction}

Tomographic imaging is based on recording projection images of an object along several directions of view. The resulting data can be interpreted as a collection of line integrals of an unknown attenuation coefficient function $f(x)$. In this work, we discretize the problem by approximating $f$ as a vectorized pixel image $\f\in\R^{N^2}$ and using the pencil-beam model for X-rays, so the indirect measurement is modelled by a matrix equation $\mbA \f= \m$. The inverse problem of reconstructing $\f$ from tomographic data is highly sensitive to noise and modelling errors, or in other words {\it ill-posed}.   

We focus on overcoming ill-posedness by enforcing sparsity of $\f$ in an orthonormal wavelet basis $\{\psi_{\gamma}\}_{\gamma\in\Gamma}$. 

In practice, the sparse reconstruction $\f_{\tiny{\text{S}}} \in \R^{N^2}$ is defined as the minimizer of this variational regularization functional:
\begin{equation}\label{eq:model}
\f_{\tiny{\text{S}}} = 
	\mbArgmin_{\f \in \R^{N^2}} \Bigg\{ \frac{1}{2}\| \mbA \f- \m \|_2^2 + \mu \sum_{\gamma\in\Gamma}|\< f,\psi_{\gamma} \>| \Bigg\}.
\end{equation}
The parameter $\mu$ in \eqref{eq:model} describes a trade-off between emphasizing more the data fidelity term or the regularizing penalty term. In general, the larger the noise amplitude in the data, the larger $\mu$ needs to be.

One popular method to solve problem \eqref{eq:model} is the so-called iterative soft-thresholding algorithm (ISTA). Such algorithm has been studied already in \cite{Lions1979}; the adaptation to sparsity-promoting inversion was introduced in \cite{Daubechies2004} and further developed in \cite{Loris2011}. Nevertheless, convergence rate for a constrained problem, such as  non-negativity constraints, is not taken into account in \cite{Daubechies2004,Loris2011}. However, in tomographic problems, enforcing non-negativity on the attenuation coefficients is highly desired. This is based on the physical fact that the X-ray radiation can only attenuate inside the target, not strengthen.
Thus, the problem we need to solve reads as:
\begin{equation}\label{eq:modelNN}
\f_{\tiny{\text{S}}} = 
	\mbArgmin_{\f \in \R^{N^2}, \, \f > 0} \Bigg\{ \frac{1}{2}\| \mbA \f- \m \|_2^2 + \mu \sum_{\gamma\in\Gamma}|\< f,\psi_{\gamma} \>| \Bigg\}.
\end{equation}
where the inequality $\f > 0$ is meant component-wise.  
In their seminal paper~\cite{chen2016primal}, Peijun Chen, Jianguo Huang, and Xiaoqun Zhang show that the minimizer of (\ref{eq:model}) can be computed using the primal-dual fixed point (PDFP) algorithm: 

\begin{equation}\label{PDFP}
\begin{split}
\y^{(i+1)} &= \mathbb{P}_C \Big(\f^{(i)} - \tau \nabla g(\f^{(i)}) -\lambda\mbW^T \v^{(i)}\Big)\\
\v^{(i+1)} &= \Big(I - \mathcal{T}_\mu \Big) \Big(\mbW \y^{(i+1)} + \v^{(i)} \Big)\\
\f^{(i+1)} &= \mathbb{P}_C \Big(\f^{(i)} - \tau \nabla g(\f^{(i)}) -\lambda\mbW^T \v^{(i+1)} \Big)
\end{split}
\end{equation}
where  $\tau$ and $\lambda$ are positive parameters, $g(\f) = \frac{1}{2}\| \mbA \f- \m \|_2^2$, the matrix $\mbW$ is a digital implementation of the wavelet transform and $\mathcal{T}$ is the
soft-thresholding operator defined by
\begin{equation}\label{eq:thresholding}
\mathcal{T}_\mu(c)=  
\begin{cases}
c + \frac{\mu}{2} &{\text {if }} x\leq -\frac{\mu}{2}\\
0 &{\text {if }} |x|<\frac{\mu}{2}\\
c - \frac{\mu}{2} &{\text {if }} x\geq -\frac{\mu}{2}.
\end{cases}
\end{equation}
Here $\mu > 0$ represents the thresholding parameter, while $\tau$ and $\lambda$ are parameters that needs to be suitably chosen to guarantee convergence. In detail, $0 < \lambda < 1/ \lambda_{\max}(\mbW\mbW^T)$, where $\lambda_{\max}$ denotes the maximum eigenvalue, and $0 < \tau < 2/\tau_{\text{lip}}$, being $\tau_{\text{lip}}$ the Lipschitz constant for $g(\f)$. Furthermore, in (\ref{PDFP}) the non-negative ``quadrant'' is denoted by $C = \R_+^{N^2}$ and $\mathbb{P}_C$ is the euclidian projection. In other words,  $\mathbb{P}_C$ replaces any negative elements in the input vector by zero.

Choosing the soft-thresholding parameter $\mu$ is analogous to the notoriously difficult problem of picking the optimal regularization parameter in Tikhonov regularization.  Many approaches for the regularization parameter selection have been proposed. For a selection of methods designed for total variation (TV) regularization see the following studies: \cite{rullgard2008,clason2010,dong2011,frick2012,wen2012,chen2014,toma2015,niinimaki2014multi}.  In this paper we introduce a novel automatic method for choosing $\mu$ based on a control algorithm driving the sparsity of the reconstruction to an {\it a priori} known ratio $0\leq {\mathcal C_{pr}}\leq 1$ of nonzero wavelet coefficients in $\f$. Our approach is based on the following idea: {\it in sparsity-promoting regularization, it is natural to assume that the} a priori {\it information is given as the percentage of nonzero coefficients in the unknown.} The idea of using the a priori known level of sparsity has been used previously \cite{kolehmainen2012,hamalainen2013sparse}, however the idea of using feedback control to achieve this is new. 

We think of the iteration \eqref{PDFP} as a \textit{plant} which takes the current threshold parameter $\mu^{(i)}$ as an input and returns $\mathcal{C}^{(i)}$, the level of sparsity in the iterate $\f^{(i)}$, as an output. 
Then, we apply a simple incremental feedback control to $\mu^{(i)}$. The feedback loop we propose is inspired by the proportional-integral-derivative (PID) controllers, which are widely used to control industrial processes \cite{aastrom1995pid,araki2009pid,bennett1993history}. If $\mathcal{C}_{pr}$ is the expected degree of sparsity, and $\mathcal{C}^{(i)}$ is the degree of sparsity at the current iterate, we change $\mu^{(i)}$ adaptively as follows:
\begin{equation}
\mu^{(i+1)} = \mu^{(i)} + \beta (\mathcal{C}^{(i)} - \mathcal{C}_{pr}), 
\end{equation}
where $\beta > 0$ is a parameter used to tune the controller. We propose a simple method for choosing $\beta$ based on the wavelet coefficients of the backprojection reconstruction, which is quick and easy to compute. If $\beta$ is chosen too large, then the controller results in an oscillating behavior for the sequence $(\mu^{(i)})_i$. On the other hand, if $\beta$ is chosen too small, reaching the expected sparsity level may take a long time.
Therefore we also account for an additional fine-tuning of the controller by exploiting the zero-crossings of the controller error 
$e^{(i)} = \mathcal{C}^{(i)} - \mathcal{C}_{pr}$.

We test our \textit{fully} automatic controlled wavelet domain sparsity (CWDS) method on both simulated and real tomographic data. The results suggest that the method produces robust and accurate reconstructions, when the suitable degree of sparsity is available.

CWDS has a connection to the following studies, which also use a parameter changing adaptively during the iterations: \cite{Bahraoui1994,Attouch1996A,Attouch1996B,Cabot2005,Rosasco2014,Rosasco2016,Hale2008}. However, our approach is different from all of them as it promotes an {\it a priori} known level of sparsity.
Also, this is not the first study which uses the wavelet transform as a regularization tool in limited data tomography. A non-exhaustive list includes~\cite{rantala2006wavelet,niinimaki2007bayesian,soussen2008reconstruction,klann2011wavelet,hamalainen2013sparse,klann2015wavelet,helin2013wavelet}. However, the proposed approach is different from the previous works, since it promotes a \textit{fully} automatic choice for 
the regularization parameter.

\section{Materials and Methods}

\subsection{Tomography setup}
%In x-ray tomography problem, the goal is to recover the density function of an unknown object from the measured projection data. In this problem, the object is a three-dimensional object and cone beam geometry is used for the X-ray source and detector.

%The mathematical model for this problem is given as follows. 
Consider a physical domain $\Omega \subset \R^2$ and a non-negative attenuation function $f:\Omega  \subset \R^2 \rightarrow\R_+$. As outlined in the Introduction, we represent $f$ by a matrix 
$\f = [\tilde{f}_{ij}] \in \R^{N \times N}$ that is later on intended as a vector belonging to $\R^{N^2}$, obtained by stacking the entries of the matrix column by column.  In X-ray tomography, the detector measures the incoming photons and the measurement data are collected from the intensity losses of X-rays from different directions or angles of view. After calibration, the measurements can be modeled as
\[
\int_{L_X} f(x) \, ds 
	= \sum_{i=1}^N \sum_{j=1}^N a_{ij} \tilde{f}_{ij},
\]
where $a_{ij}$ is the distance that a X-ray line $L_X$ travels through 
the pixel $(i,j)$.
This results in the following matrix model:
\begin{equation}
\m = \mbA \f,
\end{equation}
where the measurement matrix $\mathbf{A}=[a_{ij}] \in \R^{P \times N^2}$ contains the information about the measurement geometry, and $\m \in \R^{P}$ is the vector representing the measured data (also called \textit{sinogram}), $P$ being the number of angles of view multiplied by the number of detector cells.  

Notice that, in the following, we assume both the measurement matrix $\mbA$ and the measured data $\m$ to be normalized by the norm $\|\mbA\|$ of the matrix $\mbA$.

\subsection{2D Haar wavelets}\label{sec:Haar}
For the readers' sake of convenience, we briefly recall here the main ideas about Haar wavelets.

Consider the two real-valued functions $\varphi(x)$ and $\psi(x)$ defined on the interval $[0,1]$. Generally, $\varphi(x)$ is referred to as \textit{scaling function} and $\psi(x)$ as \textit{mother wavelet}. They are defined as follows:
\[
 \varphi(x) \equiv 1,
 \qquad 
 \psi(x) = \left\{\begin{array}{rl}
 1 &  \; \text{if} \quad 0\leq x< 1/2,\\
 -1 & \; \text{if} \quad 1/2 \leq x \leq 1.
 \end{array}\right.
\]
A \textit{Haar wavelet system} is built by appropriately scaling and translating the mother wavelet $\psi(x)$:
\[
  \psi_{jk}(x) := 2^{j/2} \psi(2^jx-k) \qquad \mbox{ for } j\leq 0, \, 0\leq k\leq 2^j-1,
\]
and the scaling function $\varphi(x)$:
\[
  \varphi_{jk}(x) := 2^{j/2} \varphi(2^jx-k) \qquad \mbox{ for } j\leq 0, \, 0\leq k\leq 2^j-1,
\]
where $\varphi(x)=0$ for $x<0$ and $x>1$. Here, $j,\, k 
\in \Z$.

It is well known that the above 1D construction leads to an orthonormal system. In 2D, we consider the standard tensor-product extension of the 1D Haar wavelet transform. In detail, a 2D Haar system is spanned by four types of functions. Three of these types have the following form: 
\begin{equation}\label{eq:threetypes}
\varphi_{jk}(x)\psi_{jk}(y),\quad 
\psi_{jk}(x)\varphi_{jk}(y), \quad
\psi_{jk}(x)\psi_{jk}(y),
\end{equation}
and the fourth type is given by $\varphi_{j_0k}(x)\varphi_{j_0k}(y)$. Notice that the fourth type describes  the coarsest scale $j_0$. The associated \textit{wavelet transform} of a function $f$ is given by
\begin{equation}\label{eq:WavTransf}
f \quad \lra \quad \mathcal{W}f = \< f,\psi_{\gamma}\>,
\qquad \gamma\in\Gamma
\end{equation}
where $\< f,\psi_{\gamma}\>$ denotes the so-called \textit{wavelet coefficients}.
Here, for notational convenience, we use the index  $\gamma\in\Gamma$ to combine together the three types \eqref{eq:threetypes} at several scales $j$ and locations $k$, and the fourth type at several locations. 

In the following, we are interested in the digital setting, \ie{}, we consider the matrix underlying the wavelet transform, which we shall denote by $\mbW \in \R^{N^2 \times N^2}$. If $\f \in \R^{N^2}$, the vector collecting all the wavelet coefficients is given by:
\begin{equation}\label{eq:WavCoeff}
   \mbW \f \; \in \, \R^{N^2},
\end{equation}
where it is clear that the matrix product $\mbW \f$ is the digital counterpart of \eqref{eq:WavTransf}.
With the above notation, the minimization problem \eqref{eq:modelNN} reads as
\begin{equation}\label{eq:modelMatrixForm}
\f_{\tiny{\text{S}}} = 
	\mbArgmin_{\f \in \R_+^{N^2}} \Bigg\{ \frac{1}{2}\| \mbA \f- \m \|_2^2 + \mu \| \mbW \f \|_1\Bigg\}.
\end{equation}

One of the main benefit of wavelets is that the transform coefficients are easy to compute and many fast algorithmic implementation are available.

For more information about the Haar wavelet transform, and its implementation, we refer to the classic text~\cite{Daubechies1992}.

\subsection{Sparsity promoting-regularization}

We consider the functional in \eqref{eq:modelNN} with 
$\{\psi_{\gamma}\}_{\gamma\in\Gamma}$ being the Haar wavelet basis as described in Subsection~\ref{sec:Haar}.
To solve the  minimization problem \eqref{eq:modelNN}, we implement the PDFP algorithm \eqref{PDFP}.

\subsection{Sparsity selection}

We assume that we have available an object $\f_{pr}$ similar to the one we are imaging. 

Given $\kappa\ge 0$, for a vector $w\in\R^{N^2}$ we define the number of elements larger than $\kappa$ in absolute value as follows: 
$$
\#_\kappa w := \#\{\,i\,\,| 1\leq i\leq N^2,\  |w_i|>\kappa\}.
$$
Now, the prior sparsity level is defined by
\[
{\mathcal{C}_{pr}} =  \frac{\#_\kappa\{\mbW \f_{pr}\}}{N^2},
\]
where $N^2$ is the total number of coefficients. In practical computations the value of $\kappa$ is set to be small but positive.

\subsection{Automatic selection of the soft-thresholding parameter $\mu$}
\label{sec:AutSel}

Assume that we know \textit{a priori} the expected degree of sparsity $0<\mathcal{C}_{pr}\leq 1$ in the reconstruction. We introduce a simple feedback loop to drive the soft-thresholding parameter $\mu$ to the desired ratio $\mathcal{C}_{pr}$ of nonzero wavelet coefficients.

The core idea is to allow $\mu=\mu^{(i)}$ to vary during the iterations by adaptively tuning it at each iteration by the following updating rule:
\[
\mu^{(i+1)} := \mu^{(i)} + \beta ({\mathcal C}^{(i)}-{\mathcal C_{pr}}),
\]
where $0\leq{\mathcal C}^{(i)}\leq 1$ is the sparsity level of the  reconstruction $\f^{(i)}$ at the $i$-th iteration. The above controller is a special case of an incremental PID-controller, where only integral control is performed. 

\subsection{The tuning parameter $\beta$}\label{TuningBeta}
Selecting the tuning parameter $\beta$ is easier than selecting the soft-thresholding parameter $\mu$. Indeed, $\beta$ has to be small enough to avoid oscillations in the sparsity $\mathcal{C}^{(i)}$ of the iterates as a function of $i$. If $\beta$ is chosen too small, this only result in a slower convergence of the algorithm.

To this purpose, we choose $\beta$ by making a suitable guess for the initial $\mu^{(0)}$. First, we compute the back-projection of the measured data to get a rough reconstruction. Back-projection is quick to compute and shows the dominant features of the target, but noise and artefacts are still predominant. As a result, the back-projection reconstruction is only good enough for estimating an initial guess for $\mu^{(0)}$, which is done by computing its wavelet coefficients. 
The initial value of the thresholding parameter $\mu^{(0)}$ is set equal to the mean of the absolute values of the $M$ smallest wavelet coefficients. In our case, we choose  
$M = n \, (1 - \mathcal{C}_{pr})$, where $n$ is the total number of wavelet coefficients. Lastly, the tuning parameter is set to be 
$\beta = \omega \mu^{(0)}$, where $\omega$ is a positive parameter. To start with a small value of $\beta$, $\omega$ is required to be small, and vice versa. 

In addition, the controller is fine tuned by detecting when the sign of difference $e^{(i)} = \mathcal{C}^{(i)} - \mathcal{C}_{pr}$ 
changes. When this happens, $\beta$ is updated by 
$\beta |e^{(i)} - e^{(i-1)}|$.
The underlying idea is that, if the desired sparsity level is crossed, that is, $e$ changes sign, either $\beta$ is far too large and oscillations have emerged, or we are already reasonably close to the optimal $\mu$ and $\beta$  can be decreased without affecting the performance too much.

%Later we call the full reconstruction image as a ground truth and denote it as $G$.

% THIS IS ONLY RELATED TO THE TV CASE WHICH IS NOT INCLUDED HERE
% Let us denote the sparsity level ${\mathcal C}$ as the number of nonzero pixel values in the images. The {\it a priori} value for the sparsity level ${\mathcal C}$ was estimated using the three photographs of the walnut as it can be seen in Figure~\ref{fig:WalnutPhoto}. The images have a size $328 \times 328$. The number of nonzero pixel values of each image was calculated. See Table~\ref{Walnut sparsity}. The percentage of the average of the three sparsity values is set as {\it a priori} sparsity level
% \begin{equation*}
% C_{pr} = 59\%. 
% \end{equation*}
%>

\subsection{Pseudo-algorithm}\label{sec:PseudoAlgorithm}
A step-by-step description of the proposed CWDS 
algorithm is summarized in Algorithm \ref{alg:Algorithm1}.

\begin{algorithm}
\caption{Controlled Wavelet Domain Sparsity Algorithm}
\begin{algorithmic}[1]
\State Inputs: measurement data vector $\mathbf{m}$, system matrix $\mathbf{A}$, parameters $\tau, \, \lambda > 0$ to ensure convergence,
{\it a priori} degree of sparsity ${\mathcal C}_{pr}$,
initial thresholding parameter $\mu^{(0)}$,
the maximum number of iterations $I_{\max}>0$, 
tolerances $\epsilon_1,\epsilon_2>0$ for the stopping rule and
control stepsize $\beta>0$.

\State $\f^{(0)} = \mathbf{0}$, $i=0$, $e=1$, and $\mathcal{C}^{(0)} = 1$ 
\While{$i<I_{\max}$ and $|e|\geq\epsilon_1$ or $d \geq \epsilon_2$}
\State $e=\mathcal{C}^{(i)}-\mathcal{C}_{pr}$
\If {$\sign(e^{(i+1)}) \neq \sign(e^{(i)})$}
\State $\beta = \beta (1 - |e^{(i+1)} - e^{(i)}|)$ 
\EndIf
\State $\mu^{(i+1)} = \max\{0, \mu^{(i)} + \beta e \}$
\State $\y^{(i+1)} = \max \{0,\f^{(i)} - \gamma \nabla g_1(\f^{(i)}) -\lambda\mbW^T \v^{(i)}\}$
% \State $\mathbf y^{i+1} = \prox_C (\mathbf f^i - \gamma \nabla g_1(\mathbf f^i) -\lambda\mcW^T \mathscr{C}^i)$
% \State $\mathscr{C}^{i+1} = (I - \prox_{\frac{\gamma}{\lambda} g_2})(\mcW \mathbf y^{i+1} + \mathscr{C}^i)$
\State $\v^{(i+1)} = (I - \mathcal{T}_{\mu^{(i)}})(\mbW \y^{(i+1)} + \v^{(i)})$
\State $\f^{(i+1)} = \max \{0,\f^{(i)}{-}\gamma \nabla g(\f^{(i)}){-}\lambda\mbW^T \v^{(i+1)}\}$
%\State $\mathbf f^{(i+1)} = \prox_C (\mathbf f^i - \gamma \nabla g(\mathbf f^i) -\lambda\mcW^T \mathscr{C}^{i+1})$
\State ${\mathcal C^{(i+1)}} = N^{-2}\#_\kappa(\mbW \f^{(i+1)})$
\State $d = \Vert \f^{(i+1)} - \f^{(i)} \Vert_2 / \Vert \f^{(i+1)} \Vert_2$
\State $i := i+1$
\EndWhile
\end{algorithmic}
\label{alg:Algorithm1}
\end{algorithm}

\section{Data Acquisition}\label{sec:DataAcquisition}
In this paper, we consider both simulated data (see Section \ref{sec:SimData}) and real data (see Section \ref{sec:RealData}).

\subsection{Simulated data}\label{sec:SimData}
\begin{figure}[h!]
\centering
\includegraphics[width=3.25cm]{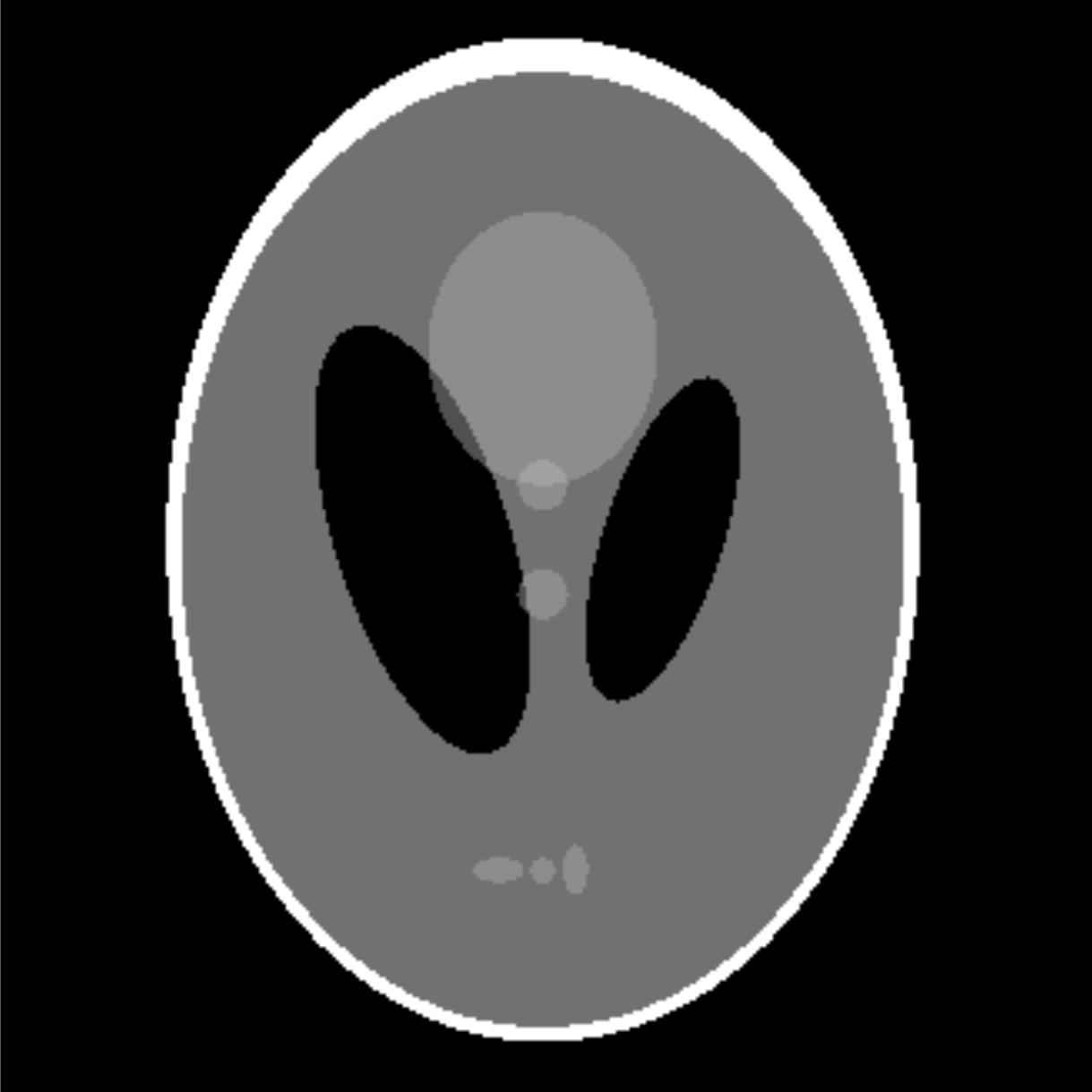}
\caption{The Shepp-Logan phantom, sized $328 \times 328$, generated with Matlab.}
\label{fig:SheppLogan}
\end{figure}

We use the Shepp-Logan phantom, available, for instance, in the Matlab Image Processing toolbox (see Figure \ref{fig:SheppLogan}). The phantom is sized $N \times N$, with $N = 328$. The projection data (\ie{}, sinogram) of the simulated phantom is corrupted by a white Gaussian process with zero mean and $0.1\%$ variance.

\subsection{Real data}\label{sec:RealData}
We use the tomographic X-ray real data of a walnut, consisting of a 2D cross-section of a real 3D walnut measured with a custom-built CT device available at the University of Helsinki (Finland). The dataset is available and freely downloadable at \url{http://fips.fi/dataset.php}. For a detailed documentation of the acquiring setup, see \cite{Hamalainen2015}.
Here we only mention that the sinogram is sized $328 \times 120$. Sinograms with different resolutions for the angle of view can be obtained by further downsampling.

\section{Numerical Experiments}
In this Section, we present preliminary numerical results in the framework of 2D fan-beam geometry.

\subsection{Algorithm parameters}
In all the experiments, we set $\lambda = 0.99$ (being $\lambda_{\max}(\mbW\mbW^T) = 1$) and $\tau = 1$ to ensure convergence. Also, we choose $\epsilon_1 = 5\times10^{-4}$ and $\epsilon_2 = 5\times10^{-4}$ for the stopping rule, and $I_{\max} = 1500$ as a safeguard maximum number of iterations (which is never attained in the results reported in Section~\ref{sec:RecRes}),  $\beta=\omega\mu^{(0)}$, where $\omega = 1$ and the values of $\mu^{(0)}$ for each experiments are shown in Table~\ref{mu0}.

\begin{table}
\caption {Initial values $\mu^{(0)}$ of the sof-thresholding parameter.}\label{mu0}
\begin{center} 
\begin{tabular}{ c | c | c | c | c | } \cline{2-4} 
   & &$120$ projections  & $30$ projections \\ \cline{2-4} 
   \multirow{2}{*}&
  Shepp-Logan & $0.0202$& $0.0195$ \\
  & Walnut &$0.0019$ & $0.0021$ \\ \cline{2-4}
\end{tabular} 
\end{center}
\end{table}
 
All the algorithms were implemented in Matlab 8.5 (R2015a) and performed on Intel Core i5 at 2.9 GHz and CPU 8GB 1867 MHz DDR3 memory. The Haar matrix $\mbW$ is generated by using Spot--A Linear-Operator Toolbox~\cite{SPOT}. The number of scales for the wavelet transform is set equal to 3 (see Figure~\ref{fig:WalnutPhoto2}).

\subsection{{\it A priori} sparsity level}
\begin{figure}[h]
\centering
\begin{tabular}{c@{\qquad}c}
\includegraphics[width=3.25cm]{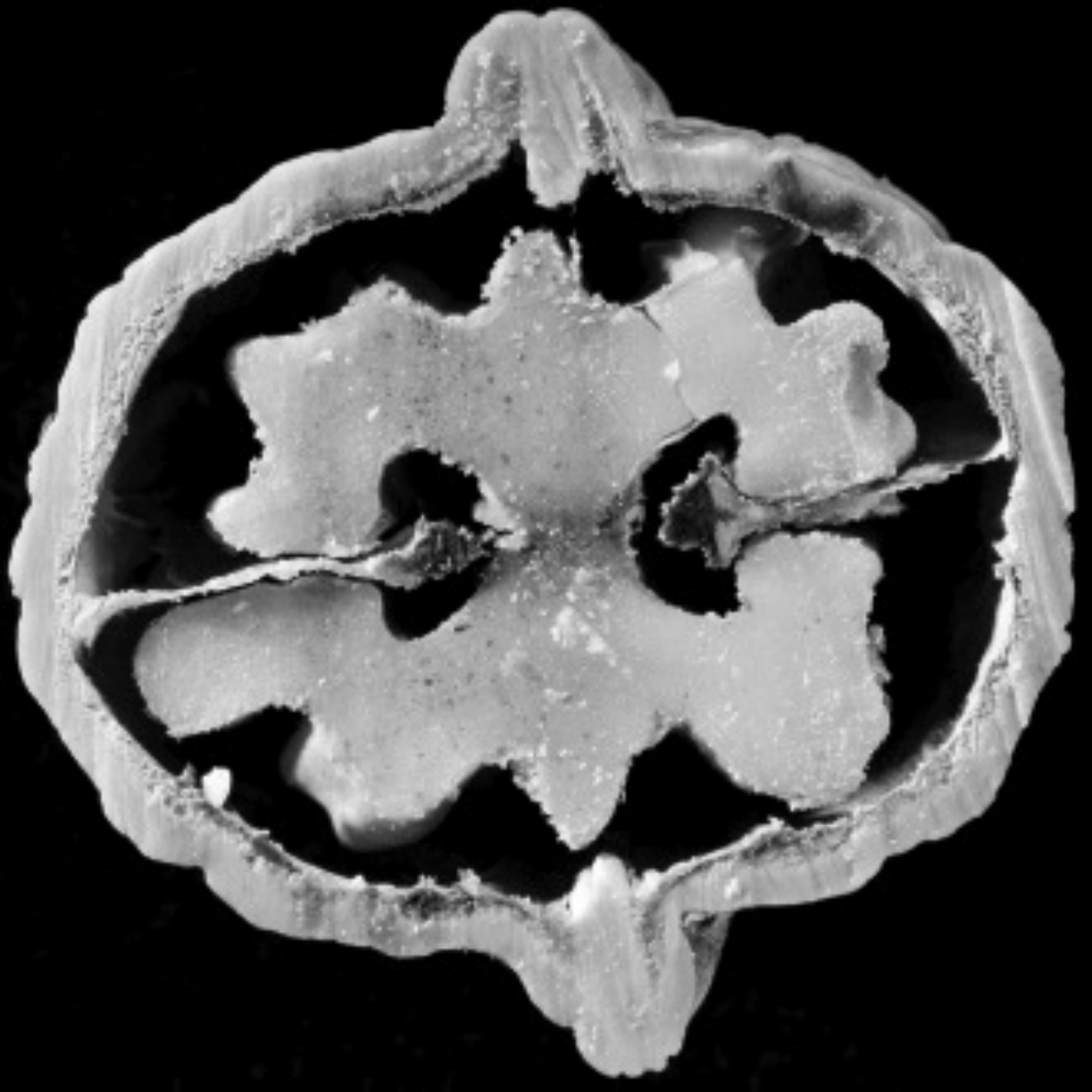}
&\includegraphics[width=3.25cm]{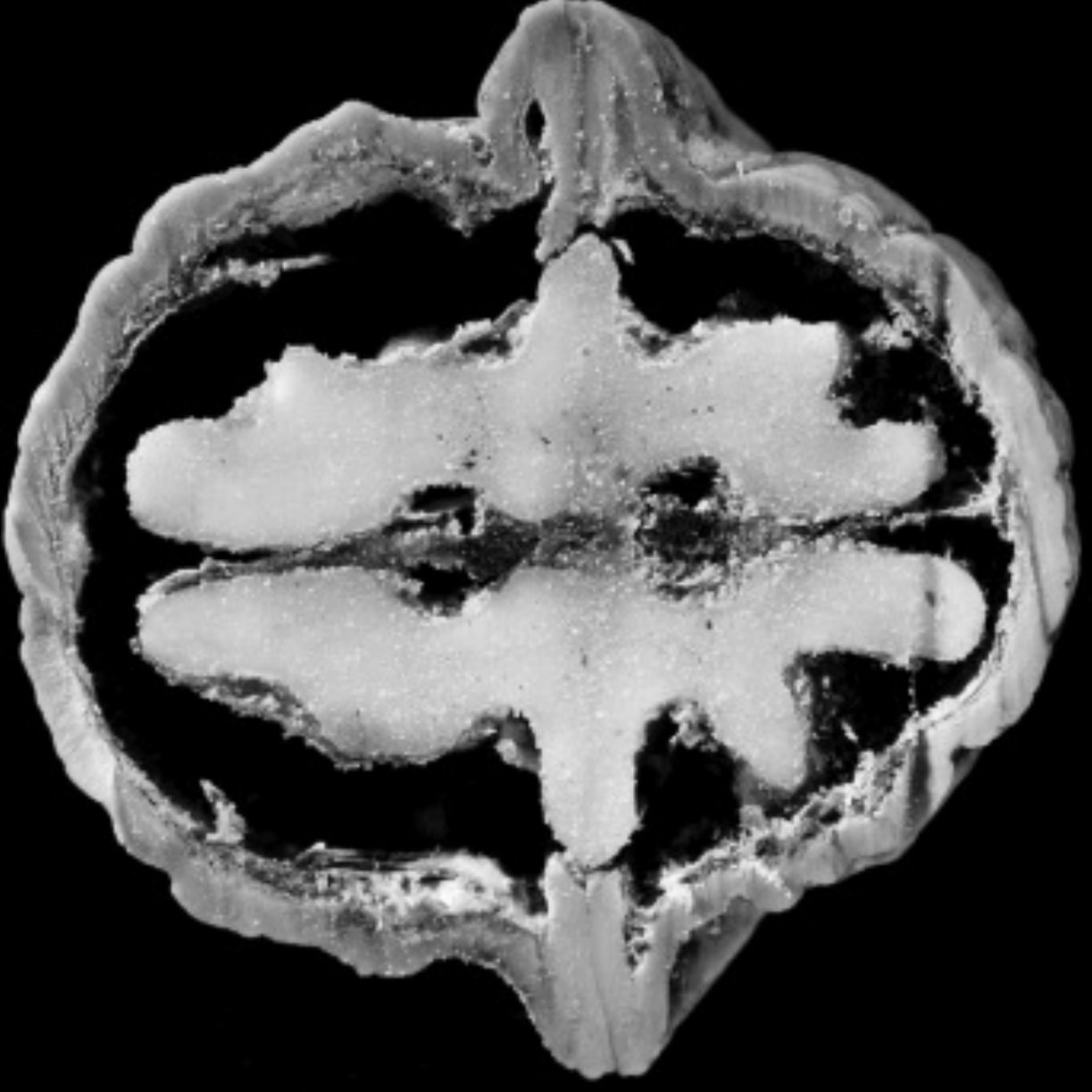}
\end{tabular}
\caption{Photographs of  walnuts split in half. The sparsity level of each image was calculated to provide the {\it a priori} information of the sparsity level for the measured walnut. The above photographs do not include the measured walnut.}
\label{fig:WalnutPhoto}
\end{figure}

\begin{figure}[h]
\centering
\includegraphics[width=7.0cm]{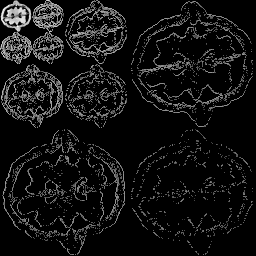}
\caption{Wavelet transform of the left photograph in Figure~\ref{fig:WalnutPhoto}. The original image is high-pass filtered, yielding the three large images. It is then low-pass filtered and downscaled, yielding an approximation image; this image is high-pass filtered to produce the three smaller detail images, and low-pass filtered to produce the final approximation image in the upper-left corner.}
\label{fig:WalnutPhoto2}
\end{figure}

To compute the desired sparsity level, we choose 
$\kappa = 10^{-6}$ for both the Shepp-Logan phantom and the walnut, and we apply the strategy outlined in Section~\ref{sec:AutSel}. In particular, for the walnut case, since we do not have at disposal the ``original'' target, we compute the sparsity level from the photographs of two walnuts cut in half (see Figure~\ref{fig:WalnutPhoto}). The \textit{a priori} sparsity level $\mathcal{C}_{pr}$ for the walnut is the average of those two sparsity levels.

For the Shepp-Logan phantom, the percentage of nonzero coefficients was estimated to be $12\%$. The percentage of the nonzero coefficients for the walnut case was estimated to be $32\%$. 

\subsection{Reconstruction results}\label{sec:RecRes}
In this Section, we present numerical results for the CWDS method, using both 
simulated and real data. As a benchmark comparison, filtered back-projection (FBP) reconstructions were also computed. For both simulated and real data, we computed reconstructions for two different resolutions of the angle of view, namely 120 and 30 projection directions, respectively.

\begin{figure}[h]
\centering
\begin{tabular}{c@{\qquad}c}
\includegraphics[width=3.25cm]{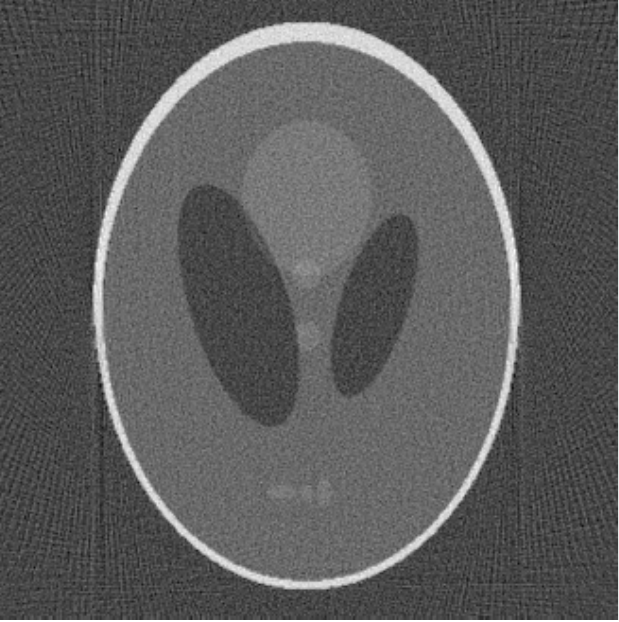}
& \includegraphics[width=3.25cm]{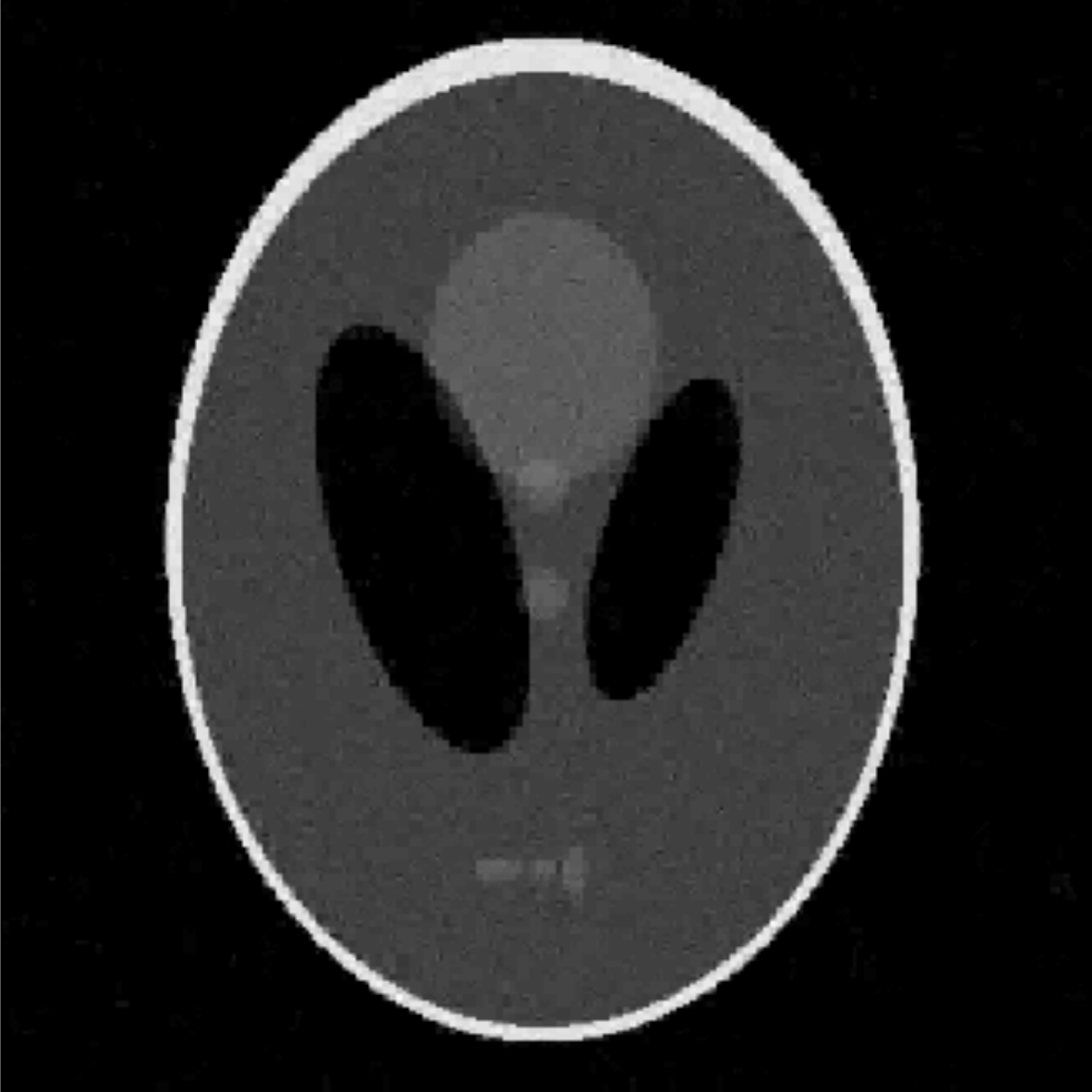} \\
(a) & (b) \\[0.5em]
\includegraphics[width=3.25cm]{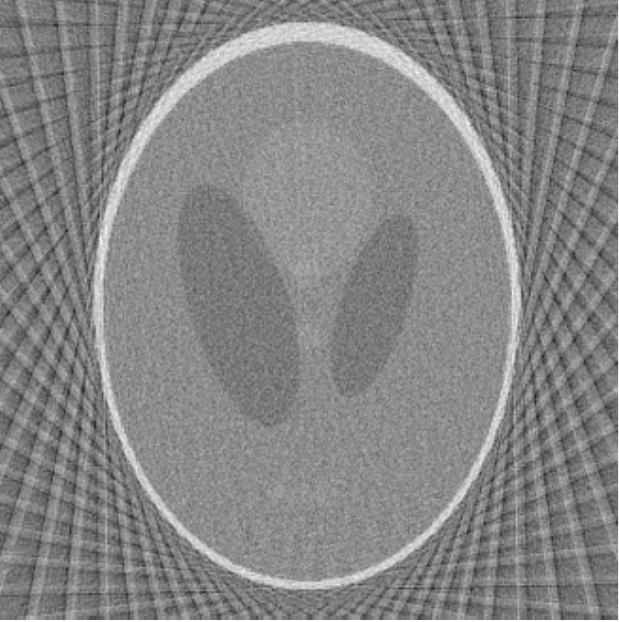}
& \includegraphics[width=3.25cm]{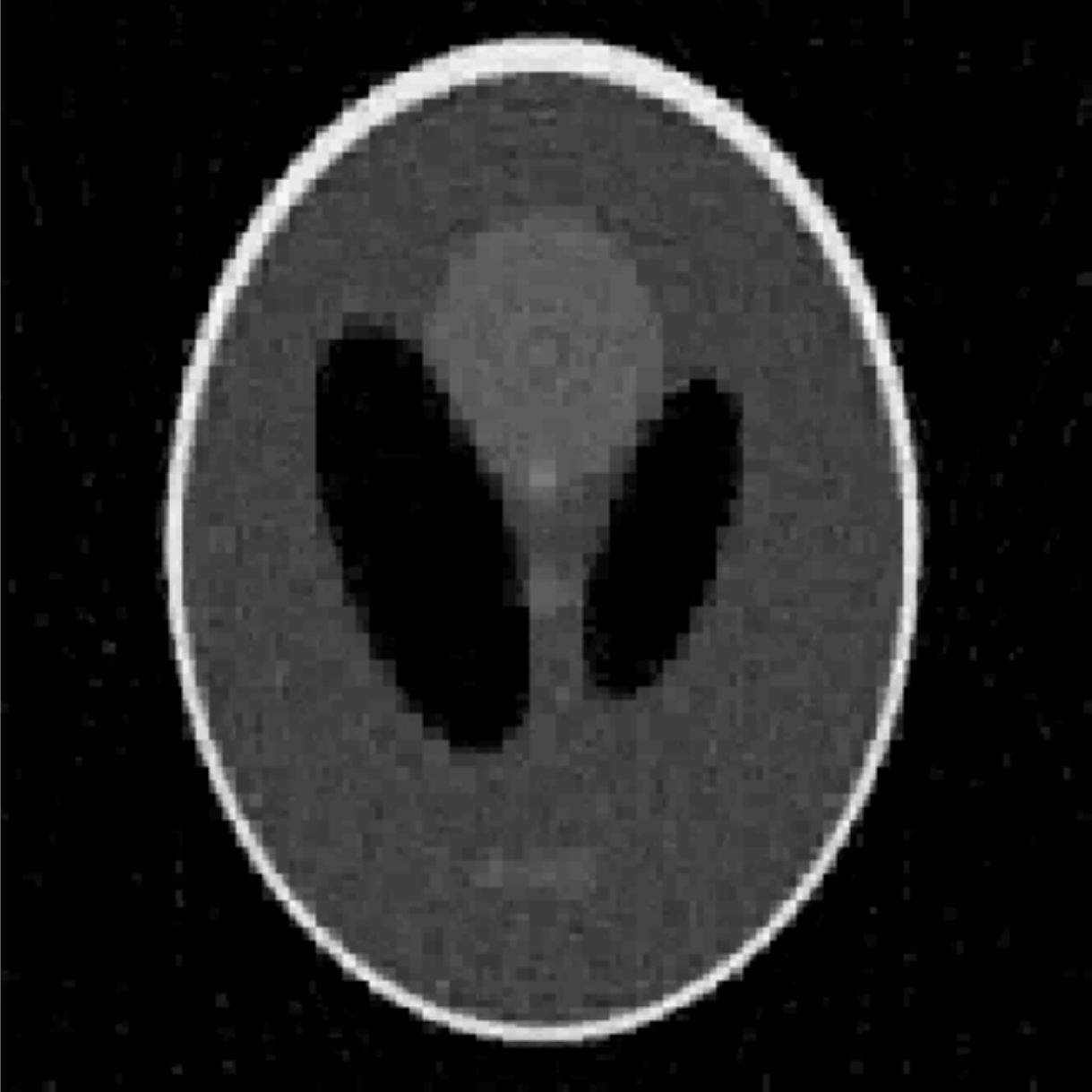} \\
(c) & (d) 
\end{tabular}
\caption{Reconstructions of the Shepp-Logan phantom using FBP with (a) 120 projections, and (c) 30 projections. Reconstructions using the wavelet based method with (b) 120 projections, and (d) 30 projections.}
\label{fig:SheppRec}
\end{figure}

\begin{figure}
\begin{tabular}{c}
\includegraphics[width=8cm]{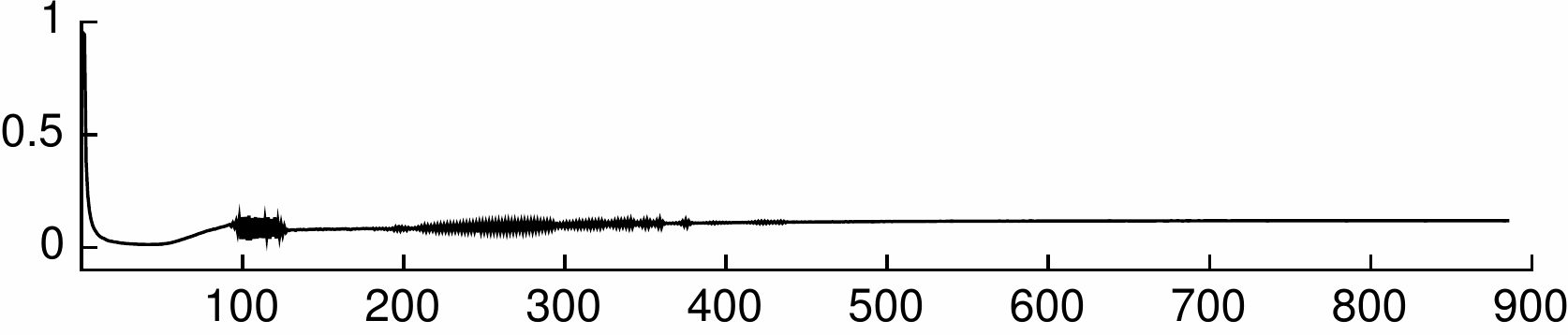} \\
\includegraphics[width=8cm]{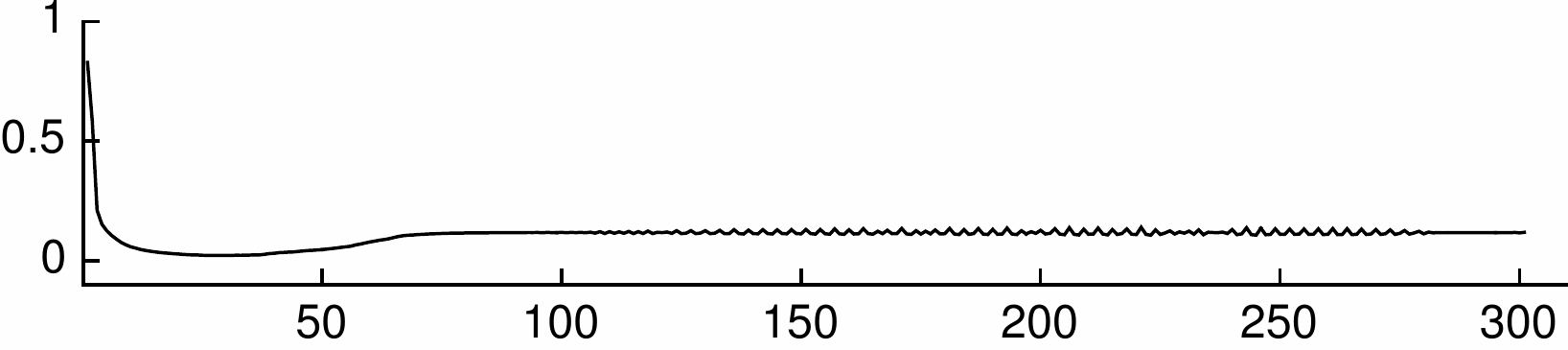}
\end{tabular}
\caption{The ratio of nonzero wavelet coefficients as the iteration progresses, for the Shepp-Logan phantom. Top: 120 projections. Bottom: 30 projections.} 
\label{fig:sparsitySheppLogan}
\end{figure}

The reconstructions of the Shepp-Logan phantom are shown in 
Figure~\ref{fig:SheppRec}. Plots of the sparsity levels, as the iteration progresses, are reported in Figure~\ref{fig:sparsitySheppLogan}. For the 120 projections case, the proposed approach converges in 885 iterations, while, in the 30 projections case, it converges in 301 iterations.  As figure of merit, we use the relative error: the obtained values are summarized in 
Table~\ref{RelativeErrorsFBP}, where we also report the values of the relative error obtained for the FBP reconstructions.

\begin{figure}[h]
\centering
\begin{tabular}{c@{\qquad}c}
\includegraphics[width=3.25cm]{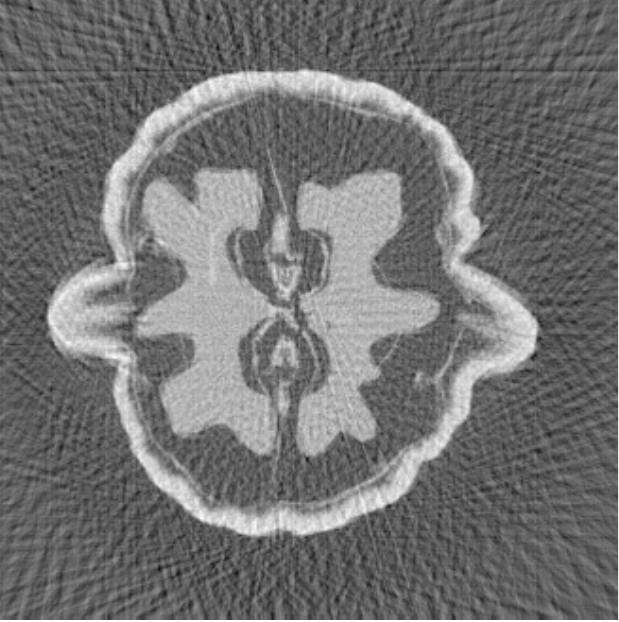} 
& \includegraphics[width=3.25cm]{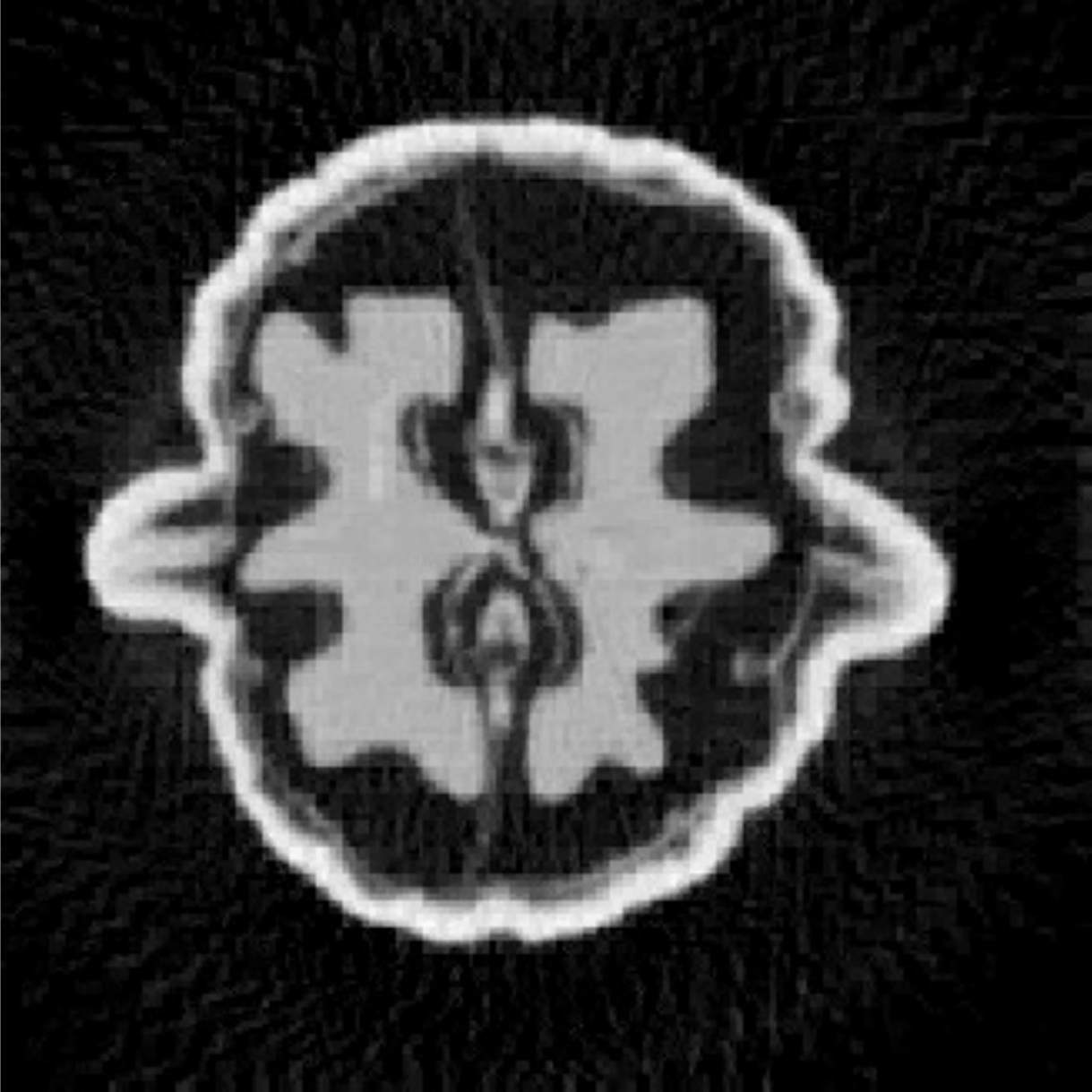} \\
(a) & (b) \\[0.5em]
\includegraphics[width=3.25cm]{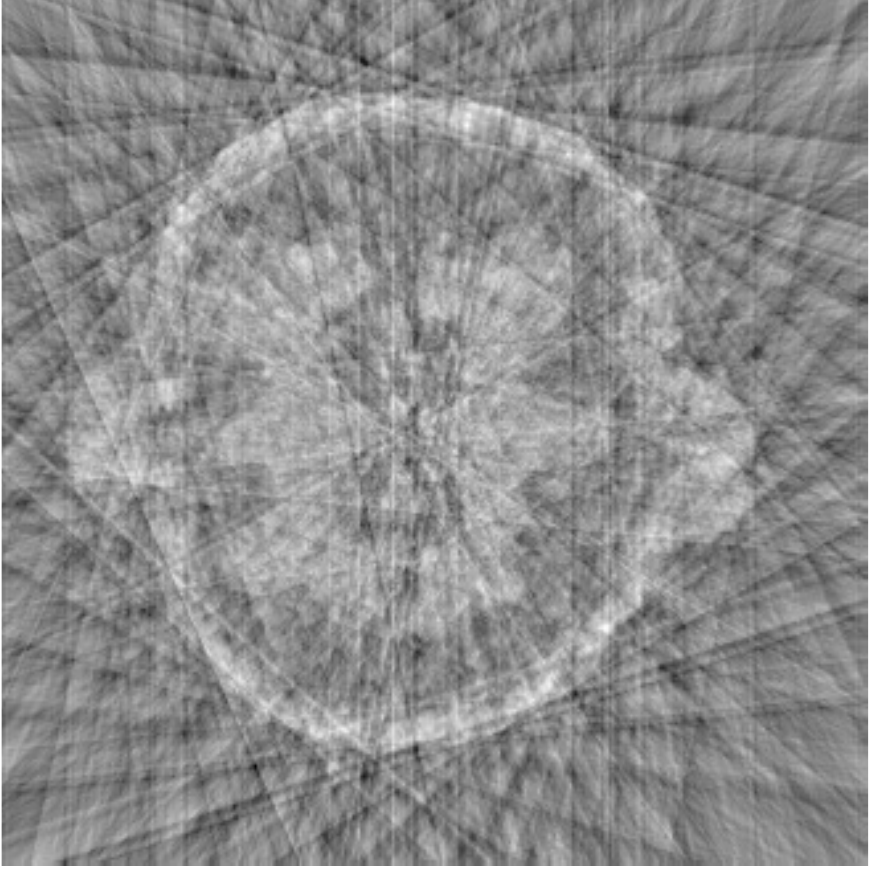}
& \includegraphics[width=3.25cm]{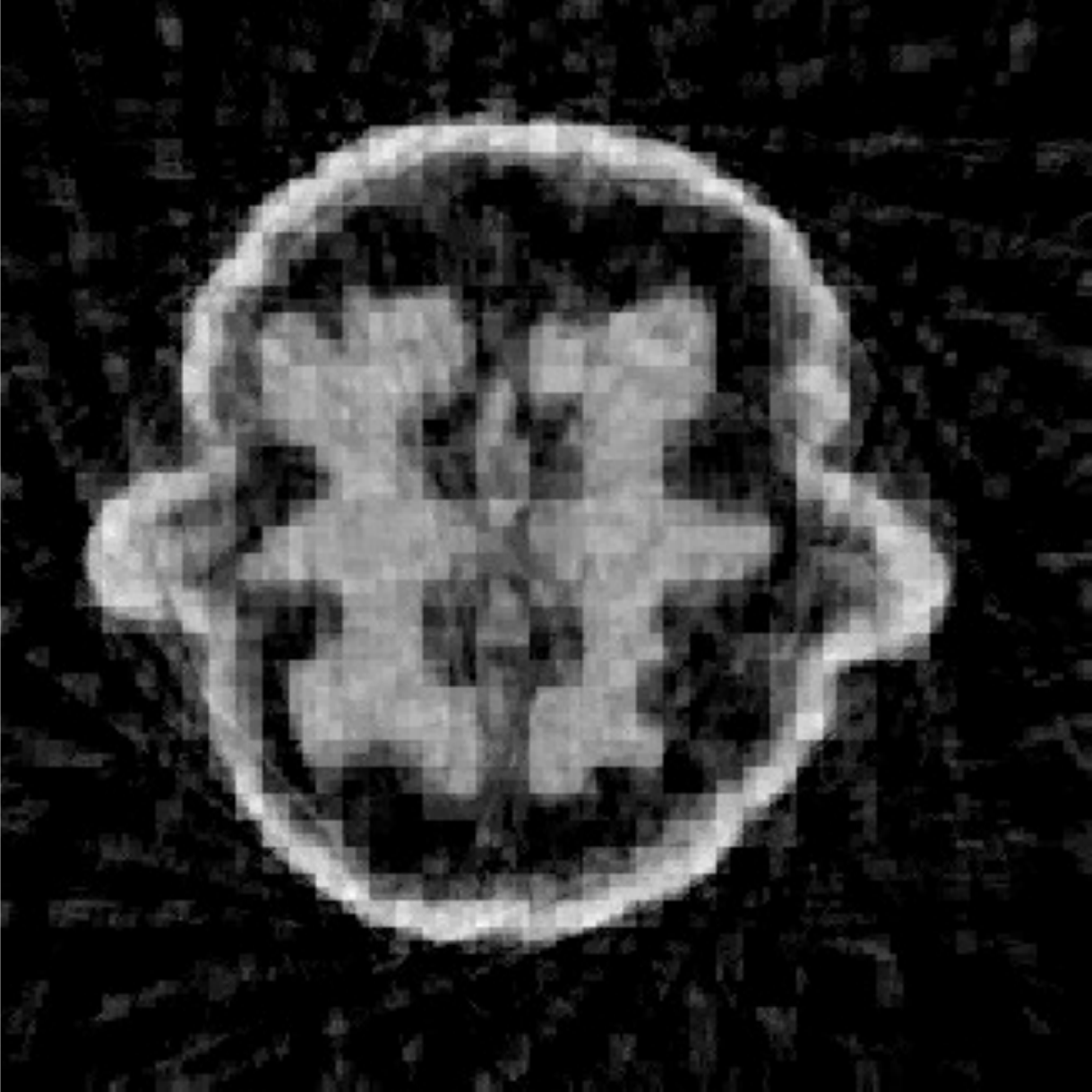} \\
(c) & (d) 
\end{tabular}
\caption{Reconstructions of the walnut using FBP with (a) 120 projections, and (c) 30 projections. Reconstructions of the walnut using the wavelet-based method with (b) 120 projections, and (d) 30 projections.}
\label{fig:WalnutRec}
\end{figure}

\begin{figure}
\begin{tabular}{c}
\includegraphics[width=8cm]{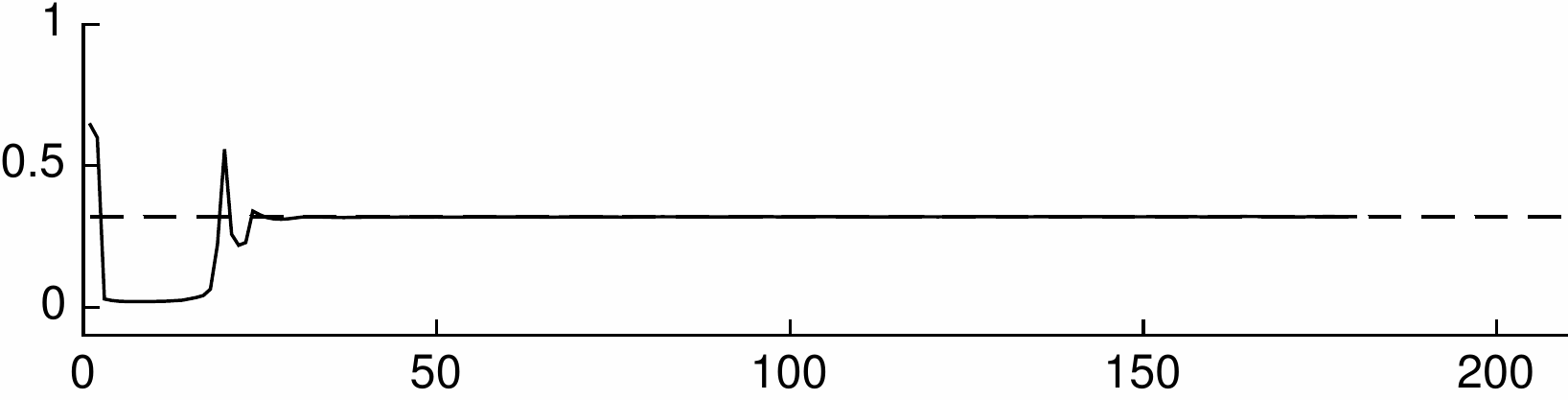} \\
\includegraphics[width=8cm]{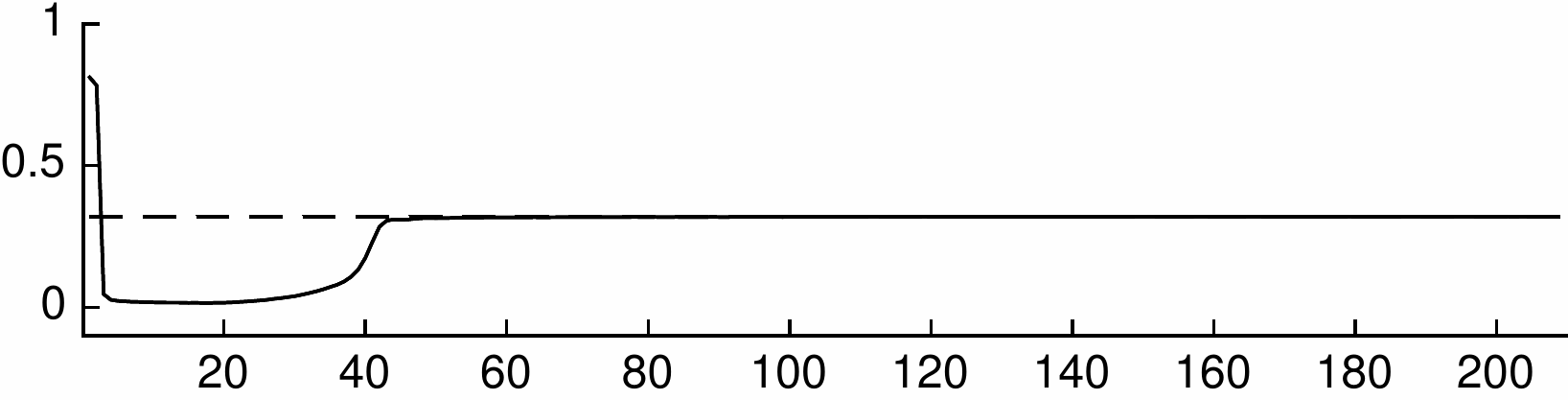}
\end{tabular}
\caption{The ratio of nonzero wavelet coefficients as the iteration progresses for the walnut case. Top: 120 projections. Bottom: 30 projections. The dashed line shows the sparsity prior $\mathcal{C}_{pr}$.} 
\label{fig:sparsityWalnut}
\end{figure}

The reconstructions for the walnut dataset, for both 120 and 30 projections, are collected in Figure~\ref{fig:WalnutRec}. The corresponding sparsity plots are shown in 
Figure~\ref{fig:sparsityWalnut}. Concerning the number of iterations to convergence, the 120 projections case required 
180 iterations, while in the 30 projections case convergence was reached in 206 iterations. 

Lastly, the computation times for all the reconstructions are reported in Table~\ref{computime}.

\begin{figure}
\begin{tabular}{c}
\includegraphics[width=8cm]{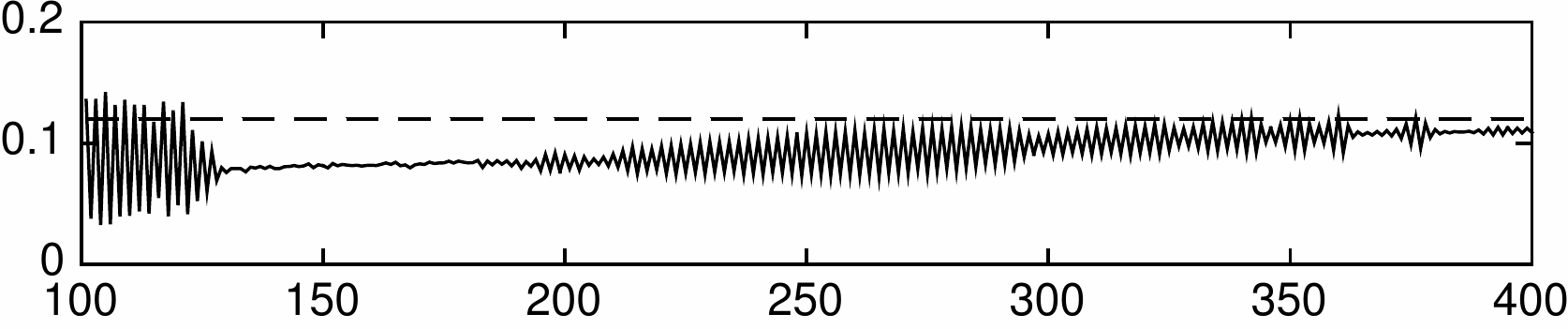} \\
\includegraphics[width=8cm]{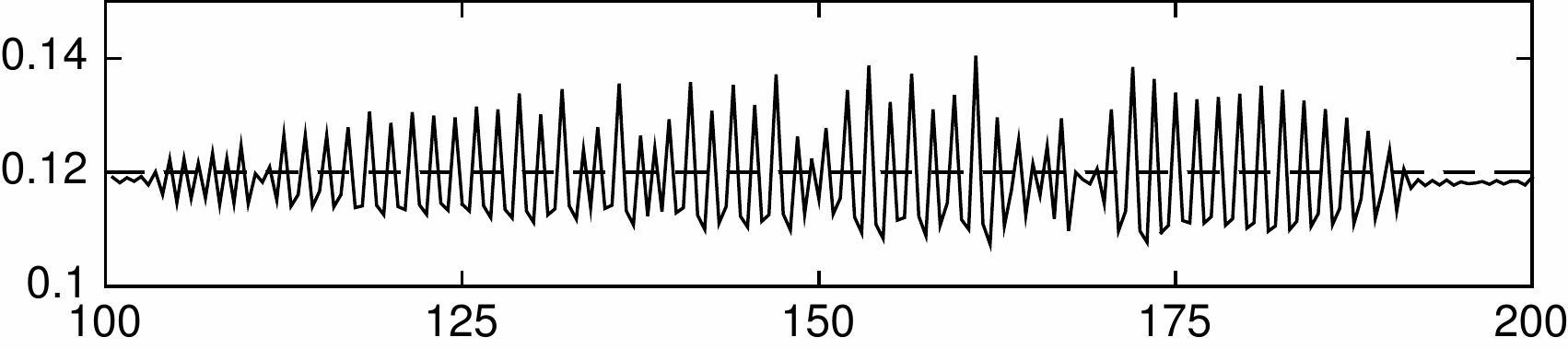}
\end{tabular}
\caption{A closer look to the sparsity level that shows oscillations for the Shepp-Logan phantom. Top: 120 projections (iterations 100 - 400). Bottom: 30 projections (iterations 100 - 200). The dashed line shows the sparsity prior $\mathcal{C}_{pr}$.} 
\label{fig:sparsityZoom}
\end{figure}

\begin{table}
\caption {The relative error of the Shepp-Logan phantom reconstructions for FBP and CWDS.} \label{RelativeErrorsFBP}
 \begin{center} 
 \begin{tabular}{ | c | c | c | c | } \hline 
   & $120$ projections  & $30$ projections  \\ \hline 
  FBP & $0.15$& $0.27$ \\
  CWDS &$0.04$ & $0.08$ \\
 \hline
 \end{tabular} 
 \end{center}
 \end{table}
 
 \begin{table}
\caption {Computation times for FBP reconstructions and CWDS reconstructions in seconds.} \label{computime}
 \begin{center} 
 \begin{tabular}{ c | c | c | c | c | } \cline{2-4} 
   & &$120$  & $30$  \\ \cline{2-4} 
   \multirow{2}{*}{walnut}&
  FBP & $0.45$& $0.09$ \\
  & CWDS &$17.40$ & $16.30$ \\ \cline{2-4}
  \multirow{2}{*}{Shepp-Logan}&
  FBP & $0.02$& $0.01$ \\
  & CWDS & $98.90$& $29.50$\\
 \cline{2-4}
 \end{tabular} 
 \end{center}
 \end{table}
 
\section{Discussion}
We presented results for both simulated and real X-ray data, also in the limited data case of only 30 projection views, with the \textit{fully} automatic CWDS method. As it can be seen in Figures~\ref{fig:SheppRec} and~\ref{fig:WalnutRec}, the reconstructions for both the Shepp-Logan phantom and the walnut data outperform the FBP reconstructions. For the Shepp-Logan case, this is confirmed by the relative errors reported in Table~\ref{RelativeErrorsFBP}.
In detail, the reconstructions using CWDS produce sharper images, with less artefacts. Overall, the quality of the reconstruction remains good even when the number of projections is reduced to 30, while, for the FBP reconstructions, streak artefacts overwhelms the reconstructions.
Finally, the presence of $\ell_1$-norm term combined with a sparsity transform, that produce denoising, and the non-negativity constraint (which is not enforced in the classical FBP scheme) definitively improves the reconstructions.

Concerning the behavior of the sparsity level for the walnut case, it can be seen in the first row of Figure~\ref{fig:sparsityWalnut} that the initial rapid oscillations decays fast. This is due to the role of the additional controller tuning $\beta$, as presented in Subsection~\ref{TuningBeta}.

For the Shepp-Logan case, it can be seen in Figure~\ref{fig:sparsitySheppLogan}, and with a closer look for some iterations in Figure~\ref{fig:sparsityZoom}, that 
the ratio of nonzero wavelet coefficients produces small oscillations for many iterations. In fact, this is a 
behavior that can appear with the proposed method: if the controller error $e$  changes sign but the absolute difference of the error of the two consecutive iterations is small, there is very little change in $\beta$. However, in 
the long run, the oscillations disappear as $\beta$ is slowly decreased. 

Future research could delve into alternative adaptive self-tuning controllers, such as the adaptive integral controller introduced in \cite{logemann1997adaptive}. Such controllers might improve the system response to unexpected disturbances and help with the oscillations caused by the slow decay of $\beta$ demonstrated in Figure~\ref{fig:sparsityZoom}. Additionally careful analysis of the dynamics of the algorithm \eqref{PDFP} is required to see if convergence of CWDS can always be guaranteed with the methods presented in this paper. 

Anyhow, what is remarkable is that, for all numerical experiments, the sparsity level eventually converges to the desired sparsity level $\mathcal{C}_{pr}$.

\section{Conclusion}
In this paper, we proposed a new approach in tuning the regularization parameter, in this case the sparsity level of the reconstruction in the wavelet domain. CWDS seems to be a promising strategy, especially in real life applications where the end-users could avoid manually tuning the parameters. 

In the case of sparsely collected projection data, the \textit{fully} automatic CWDS outperforms the conventional FBP algorithm in terms of image quality (measured as relative RMS error).

\section*{Acknowledgments}
This work was supported by the Academy of Finland through the Finnish Centre of Excellence in Inverse Problems Research 2012--2017 (Academy of
Finland CoE-project 284715).

% \ifCLASSOPTIONcaptionsoff
%   \newpage
% \fi
%\bibliographystyle{IEEEtran}
\bibliographystyle{unsrt}
\bibliography{Controlled Wavelet Domain Sparsity for X-ray Tomography_ArXiv}

\end{document}